\documentclass[a4paper,12pt]{article}
\usepackage{amsfonts}
\usepackage{multirow}
\usepackage{graphicx}

\begin{document}

\title{Banishing divergence Part 1: Infinite numbers as the limit of sequences of real numbers.}
\author{David Alan Paterson\\
CSIRO CMSE\\
Graham Rd, Highett, 3090\\
Australia}
\date{\today}
\maketitle

\begin{abstract}
Sequences diverge either because they head off to infinity or because they oscillate. Part 1 constructs a non-Archimedean framework of infinite numbers that is large enough to contain asymptotic limit points for non-oscillating sequences that head off to infinity. It begins by defining Archimedean classes of infinite numbers. Each class is denoted by a prototype sequence. These prototypes are used as asymptotes for determining leading term limits of sequences. By subtracting off leading term limits and repeating, limits are obtained for a subset of sequences called here ``smooth sequences". $\mathbb{I}_n$ is defined as the set of ratios of limits of smooth sequences. It is shown that $\mathbb{I}_n$ is an ordered field that includes real, infinite and infinitesimal numbers.
\end{abstract}
\section{Introduction}

Applied mathematicians deal with divergence on a daily basis. Pure mathematicians deal with infinity on a daily basis. This paper is the first step towards a mathematical utopia in which ``divergence" ceases to exist and every infinite sequence of real numbers has a limit. To arrange this requires a new type of asymptotic limit. This new type of limit requires an enumeration of Archimedean classes in a non-Archimedean framework containing finite, infinite and infinitesimal numbers.

There are currently three main independent ways to define finite, infinite and infinitesimal numbers~\cite{1}: axiomatically through the hyperreals $^*\mathbb{R}$ and nonstandard analysis~\cite{2,3,4}, from a simplification of the Dedekind cut through the surreal numbers $\bf{N_O}$~\cite{5,6,7,8}, and by consideration of the geometric continuum by Veronese and Hilbert~\cite{9,10,11}. What has been missing until now has been a definition of infinite numbers as asymptotic limits of sequences, in the same way that reals are Cauchy limits on rationals. That is a gap that Part 1 of this paper fills.

Any reader who has done mathematical work on non-Archimedean numbers will soon be asking why this paper hasn't defined the infinite numbers using power series, as has been done by Hahn, Levi-Civita, Laurent, and many others~\cite{9,18,22,Hahn,23}. The reason is simple: infinite power series produce too many numbers, and finite power series don't produce enough. Infinite power series produce numbers like $\sin(\omega)$ that are finite but not real. Real numbers are the largest finite field of ordered numbers, so any field that contains numbers that are finite but not real cannot be ordered, at least not in the usual sense. A finite power series based on a number with any particular Cantor cardinal cannot produce a number with any other Cantor cardinal and is therefore incomplete. The approach here is to use a finite series but expand the set of basis functions enormously to try to walk the tightrope between too few and too many numbers.

Part 1 of this paper does not attempt to prove any significant theorem, but instead uses pure mathematics to construct a methodology for use in Part 2~\cite{Part2}. Part 2 of this paper is applied mathematics, and deals with limits of functions that oscillate, as well as with applications of infinite sequences. As much as possible of the material presented here is new.

\section{Preliminaries}

It is well known that the ordinal numbers are non-commutative, that $\omega = 1 + \omega \neq \omega + 1$. Yet there are some papers in which the ordinal numbers are either explicitly or implicitly assumed to be commutative. These include Klaua(1994)~\cite{12} and all the papers on surreal numbers~\cite{5,6,7,8}.

The statement of this non-commutative property dates back to Cantor (1895)~\cite{13}.

$\omega$ is associated with the ordered set $\{a_1, a_2, a_3, \ldots\}$. Consider a single extra element $\{b\}$.
$1 + \omega$ comes from $\{b, a_1, a_2, a_3, \ldots\}$.
$\omega + 1$ comes from $\{a_1, a_2, a_3, \ldots , b\}$.
Cantor then boldly asserts that $\omega = 1 + \omega \neq \omega + 1$.

In the present paper the approach is different. The statement $\{b, a_1, a_2, a_3$ $, \ldots\} = \{a_1, a_2, a_3, \ldots\}$ can only be true if infinite commutativity within the sequence is allowed. But the statement $\{a_1, a_2, a_3, \ldots , b\} \neq \{a_1, a_2, a_3, \ldots\}$ can only be true if infinite commutativity within the sequence is prohibited. Let only finite commutativity within the sequence be allowed. Then $\{b, a_1, a_2, a_3, \ldots\} \neq \{a_1, a_2, a_3, \ldots\}$. Further, allow finite commutativity of the form $\{b, a_1, a_2, a_3, \ldots\} = \{b\} + \{a_1, a_2, a_3, \ldots\} = \{a_1, a_2, a_3, \ldots\} + \{b\} = \{a_1, a_2, a_3, \ldots , b\}$.

When that is done, the result is a set of \textbf{commutative ordinals}, and  $\omega \neq 1 + \omega = \omega + 1$. Throughout this paper, the ordinals are taken to be commutative. Ehrlich ascribes the introduction of commutative ordinals to Hessenberg (1906)~\cite{14}. The insistence on finite commutativity also unlinks ordered sequences from the process of bijection~\cite{bijection}, which is a good thing because the process of bijection has the ability to completely destroy the ordering of an ordered sequence. By eliminating the option of using bijection, all classical paradoxes involving infinity are immediately solved.

If you're uncomfortable with that definition of $\omega$, choose whichever you're most comfortable with:
\begin{itemize}
  \item	$\omega$ is the first ordinal number
  \item	$\omega$ is the surreal number $\{\mathbb{Z}|\}$
  \item	$\omega = S(\mathbb{N})$ where $S()$ is the successor function from ZFC theory~\cite{15}
  \item	$\omega = \{\mathbb{N}\}$, the ordered set of natural numbers
\end{itemize}
In each case, the important point is that $\omega$ is greater than every natural number.\\

Big O notation is used to describe the limiting behaviour of a function when the argument tends towards a particular value of infinity~\cite{16}. Big O notation characterizes functions according to their growth rates: different functions with the same growth rate may be represented using the same O notation. 

A formal definition is~\cite{16}:
\begin{itemize}
  \item	$f(x) = O(g(x))$ if and only if there exists a positive real number $M$ and a real number $x_0$ such that $|f(x)| \le M |g(x)|$ for all $x > x_0$. 
\end{itemize}
As originally applied to operations counts in computer algorithms, $f$ is a mapping from natural numbers to natural numbers and $g$ is a mapping from natural numbers to reals. ``Big O notation" is thus a statement about sequences. For the purposes of the present paper, the definition is modified to make it symmetric. From the above definition it is perfectly valid to write $n = O(\exp(n))$ although $\exp(n) \neq O(n)$. The following definition forces $n \neq O(\exp(n))$:
\begin{itemize}
  \item	$f(x) = O(g(x))$ if and only if there exists a positive real number $M$ and a real number $x_0$ such that $|f(x)| \le M |g(x)|$ and $|g(x)| \le M |f(x)|$ for all $x > x_0$.
\end{itemize}
Call this \textbf{commutative big O notation}. With commutative big O notation, $f = O(g)$ implies $g = O(f)$.

\section{Archimedean classes and prototypes}

Different sources use different definitions of the Archimedean axiom. These definitions tend to be centred on one of the following three types of definitions. In each, $G$ is an ordered semigroup, which includes ordered groups and fields. The three types of definitions are:
\begin{enumerate}
  \item $\forall a,b\in G \ni m\in \mathbb{N} : ma>b \wedge mb>a$
  \item $\forall 0<a\in G <b\in G \ni m\in \mathbb{N} : ma>b$
  \item If $G$ is a group, i.e. $\forall a\in G : 1/a\in G$, then the following also applies:\\
$\forall a\in G \ni m\in \mathbb{N} : |a|<m$
\end{enumerate}
The Archimedean axiom can be derived from the axiom of Dedekind completeness, sometimes called the ``Least-upper-bound property"~\cite{17}. The inverse is not true, the rational numbers obey the Archimedean axiom but are not Dedekind complete. There exists an altered version of the axiom of Dedekind completeness that does not imply the Archimedean axiom. There's a well-known theorem due to Hilbert that the largest ordered field that obeys the Archimedean axiom is isomorphic to the real numbers.~\cite{2}

Any ordered field that does not satisfy the Archimedean axiom is called non-Archimedean~\cite{18}. Any ordered field that includes both finite and infinite numbers, or both finite and infinitesimal numbers, is non-Archimedean. Within any non-Archimedean field will be ordered collections of numbers that obey the Archimedean axiom. These are known as Archimedean classes~\cite{19}. 

Define an \textbf{Archimedean class} by the following
\begin{itemize}
  \item	Two numbers $a$ and $b$ have the same Archimedean class if there exists positive integer $m$ such that $m|a|>b$ and $m|b|>a$.
\end{itemize}
Clearly, zero can't be in an Archimedean class with any other number. This immediately implies that the reals, rationals and integers are not Archimedean classes. But $\mathbb{R}\backslash 0$, $\mathbb{Q}\backslash 0$ and $\mathbb{Z}\backslash 0$ are Archimedean classes. This definition allows each Archimedean class to be complete under the operations of multiplication and division, which is not possible in an ordered field because division by zero is not allowed. It also means that Archimedean classes are not complete under addition and subtraction.\\

The use of commutative big O notation to enumerate Archimedean classes of a non-Archimedean ordered field may be new.
 
\emph{Lemma}. $f=O(g)$ if and only if $f$ and $g$ are members of the same Archimedean class.
The proof is trivial.

Typical examples of Big-O notation are $O(0)$, $O(1/n)$, $O(1)$, $O(\log n)$, $O(n^{1/3})$, $O(n)$, $O(n^{3.5})$, $O(2^n)$, $O(n^n)$. Each of these examples defines a different Archimedean class. The real numbers are all members of $O(1)$ or $O(0)$. These examples show the way to creating a very large number of infinite and infinitesimal Archimedean classes.\\

Let each class be denoted by a typical member, called a \textbf{prototype} $p_n$, and this is used to identify the whole class, thus negating the need to repeat $O()$ every time. Then take the limit as $n\to\omega$ and call the result $p$. The following defines a set of Archimedean classes:\\

Part 1. The infinite Archimedean classes\\

The infinite Archimedean classes are generated by $\omega$ and by closure under a finite number of operations of:
\begin{itemize}
\item $\ln( )$, $\exp( )$, finite positive power, a multiplication of two unequal classes, and a division of a larger class by a smaller; limiting the use of ln to only $\ln(\omega)$ and $\ln(\ln(\ldots))$.
\item $\exp\exp\exp\ldots() = \exp^{\omega}()$ and $\ln\ln\ln\ldots() = \ln^{\omega}()$. Call this the \textbf{cardinal jump rule}.
\item A feedback rule that can't be defined until later when the infinite limit ``lim" has been defined. Under certain conditions a class to a power of ``1/lim" can be a class.
\end{itemize}

Part 2. Other Archimedean classes
\begin{itemize}
\item The infinitesimal Archimedean classes are 1/the-infinite-Archimedean-classes.
\item The reals: 1.
\item Zero. This contains only one member.
\end{itemize}
The above definition specifically excludes all cases of the sum of two Archi-medean classes, and the product of an Archimedean class and a finite number. 

``Closure" means that duplicates are equivalenced.

In place of $\ln( )$ and $\exp( )$, any other logarithmic and inverse logarithmic function can be substituted in the definition, for example $\log_2( )$ and $2^{( )}$.

Set aside the cardinal jump rule for discussion at the end of this section, and set the feedback rule aside until the next section on limits.\\

\emph{Lemma}: Every prototype of an Archimedean class other than 0 is greater than zero.

\emph{Proof}: $\omega$ is greater than every natural number so $\omega > 0$. Infinite Archimedean types are derived by a finite number of iterative steps from $\omega$ so use induction. Move from one step, where $p_1$ and $p_2$ are positive, to the next step. $p = \ln(p_1)$ or $\exp(p_1)$ or $p_1^\alpha$ for $\alpha > 0$ or $p_1p_2$ or $p_1/p_2$. All the options other than $p = \ln(p_1)$ force $p > 0$. Let $p = \ln(\ldots\ln(x)\ldots)$  where there are $n$ applications of $\ln()$. Then $p > 0$ if and only if $x > y = \exp(\ldots\exp(0)\ldots)$. Because there are a finite number of applications of $\exp()$, $y$ is a finite real number. $\lceil y\rceil>y$ is a natural number so $x = \omega > \lceil y\rceil>y $ and so $p > 0$. So all infinite $p$ are positive and so all infinitesimal $1/p > 0$. The only remaining prototype is $1 > 0$.

\emph{Lemma}: Every prototype of an infinite Archimedean class is greater than one.

\emph{Proof}: $\omega$ is greater than every natural number so $\omega > 1$. Infinite Archimedean types are derived by a finite number of iterative steps from $\omega$. Use induction. Move from one step, where $p_1$ and $p_2$ are greater than 1, to the next step. $p = \ln(p_1)$ or $\exp(p_1)$ or $p_1^\alpha$ for $\alpha > 0$ or $p_1p_2$ or $p_1/p_2$ where $p_1 > p_2$. All the options other than $p = \ln(p_1)$ force $p > 1$. Let $p = \ln(\ldots\ln(x)\ldots)$ where there are $n$ applications of $\ln()$. Then $p > 1$ if and only if $x > y = \exp(\ldots\exp(1)\ldots)$. Because there are a finite number of applications of $\exp()$, $y$ is a finite real number. $\lceil y\rceil>y$ is a natural number so $x = \omega > \lceil y\rceil>y$  and so $p > 1$. So all infinite $p$ are greater than one.

\emph{Lemma}: The set of prototypes of Archimedean classes other than 0 are closed under multiplication and division.

\emph{Proof}: Let $p > q > 1$ be two infinite Archimedean classes with infinitesimal counterparts $1/p$ and $1/q$. Then $p > q > 1/q > 1/p$. By the definition of an infinite Archimedean class, $pq$ and $p/q$ are infinite Archimedean classes. Therefore so are $p(1/q)$ and $p/(1/q)$. The inverses of these must be infinitesimal Archimedean classes $1/(pq) = (1/p)(1/q) = (1/p)/q$ and $1/(p/q) = (1/p)q = (1/p)/(1/q)$. This exhausts all possibilities of pairwise multiplication and division that include both $p$ and $q$. The Archimedean prototype ``1" multiples and divides infinite and infinitesimal classes as $1p = p, 1(1/p) = 1/p, 1/p = 1/p, 1/(1/p) = p$ and it is produced by $p/p = p(1/p) = 1(1) = 1/1 = 1$. That doesn't quite exhaust all possibilities yet. $p(p) = p/(1/p) = p^\alpha$ where $\alpha = 2$ is an infinite Archimedean class by definition and its inverse $(1/p)(1/p) = (1/p)/p = 1/p^\alpha$ is an infinitesimal Archimedean class. That totals 10 each multiplications and divisions and covers all possibilities.\\

Let's look at the ordering of infinite Archimedean classes for the first few generations.\\

Generation 0.
$\omega$\\

Generation 1.
$\ln(\omega), \omega^\alpha, \exp(\omega)$
for $\alpha \in \mathbb{R}^+$.

\emph{Lemma}: $\ln(\omega) < \omega^\alpha < \exp(\omega)$ for all $\alpha \in \mathbb{R}^+$.

\emph{Proof}: For given $\alpha$, the value of $x^\alpha < \exp(x)$ when $x/\ln(x) > \alpha$. The value of $x/\ln(x)$ can be made arbitrarily large by making $x$ sufficiently large finite number. And $\omega$, being infinite, is always sufficiently large.
For given $\alpha$, the value of $\ln(x) < x^\alpha$ when $\ln(\ln(x))/\ln(x) < \alpha$. Let $y = \ln(x)$. The inequality is always satisfied for a sufficiently large finite $y$ and hence a sufficiently large finite $x = \exp(y)$. $\omega$, being infinite, is always sufficiently large.

The number $\alpha > 0$ takes the same ordering as the reals. The construction $f^\alpha$ absorbs $f$ as a special case.\\

Generation 2.
If $f$ and $g$ are two members of previous generations then, subject to restrictions, the next generation contains $\ln(g), g^\alpha, \exp(g)$ and $fg$ and $f/g$ when $f > g$. In the following three tables, all prohibited operations and duplicates have been removed.

\setlength{\tabcolsep}{1 pt}
\begin{table} [htb]
\caption{Generation 2. $f(g(\omega))$.}
\centering
\begin{tabular}{|cc|ccc|}
\hline
 & & & $f()$ &\\
 & & $\ln()$ & $()^\alpha$ &	$\exp()$\\\hline
& $\ln(\omega)$ & $\ln(\ln(\omega))$ & $\ln(\omega)^\alpha$ &\\
$g()$ & $\omega^\alpha$ & & & $\exp(\omega^\alpha)$\\
& $\exp(\omega)$ & & $\exp(\alpha\omega)$ & $\exp(\exp(\omega))$\\\hline
\end{tabular}
\end{table}
\begin{table} [htb]
\caption{Generation 2. $fg$.}
\centering
\begin{tabular}{|cc|cc|}
\hline
 & & $f(\omega)$ &\\
 & & $\omega^\alpha$ & $\exp(\omega)$ \\\hline
$g(\omega)$ & $\ln(\omega)$ & $\omega^\alpha \ln(\omega)$ & $\exp(\omega)\ln(\omega)$ \\
 & $\omega^\alpha$ & & $\exp(\omega)\omega^\alpha$ \\\hline
\end{tabular}
\end{table}

\begin{table} [htb]
\caption{Generation 2. $f/g$.}
\centering
\begin{tabular}{|cc|cc|}
\hline
 & & $f(\omega)$ &\\
 & & $\omega^\alpha$ & $\exp(\omega)$ \\\hline
$g(\omega)$ & $\ln(\omega)$ & $\omega^\alpha/\ln(\omega)$ & $\exp(\omega)/\ln(\omega)$ \\
 & $\omega^\alpha$ & & $\exp(\omega)/\omega^\alpha$ \\\hline
\end{tabular}
\end{table}

\textbf{The fg and f/g rule}. The ordering of $fg$, and $f/g$ is immediate. After all other ordering is done, first subdivide using the ordering of $f$, then subdivide using the ordering of $g$ with large $g$ to small in $f/g$, then $f$, then small $g$ to large of $fg$.

This completes the 12 entries from the combined generations 0, 1 and 2. With $\alpha = 1$ these have the following order: $\ln(\ln(\omega)) < \ln(\omega)^\alpha < \omega^\alpha/\ln(\omega) < \omega^\alpha < \omega^\alpha \ln(\omega) < \exp(\omega)/\omega^\alpha < \exp(\omega)/\ln(\omega) < \exp(\alpha\omega) = \exp(\omega^\alpha) < \exp(\omega)\ln(\omega) < \exp(\omega)\omega^\alpha < \exp(\exp(\omega))$.

The equality occurs because the two are subsets of the same collection of Archimedean classes. For instance $\omega$ is a subset of $\alpha\omega^\beta$. The ordering within each $\alpha\omega^\beta$ collection is given first by $\beta$ and second by $\alpha$. Between collections the construction $\alpha\omega^\beta$ is considered monolithic.\\

Generation 3.
The ordering of the entries from combined generations 0, 1, 2 and 3 for $\alpha,\beta,\gamma,\delta=1$, yields (after deleting duplicates and after slight simplification using collections):

$\ln(\ln(\ln(\omega)))<\ln(\ln(\omega))^{\alpha}<\ln(\omega)^{\alpha}/\ln(\ln(\omega))
<\ln(\omega)^{\alpha}<\ln(\omega)^{\alpha}\ln(\ln(\omega))$

$<\omega^\alpha/\ln(\omega)/\ln(\ln(\omega)) <\omega^\alpha/\ln(\omega)^{\beta}<\omega^\alpha\ln(\ln(\omega))/\ln(\omega)
<\omega^\alpha/\ln(\ln(\omega))$ $<\exp(\alpha\ln(\omega)^\beta)
<\omega^\alpha\ln(\ln(\omega))<\omega^\alpha\ln(\omega)/\ln(\ln(\omega))
<\omega^\alpha\ln(\omega)^{\beta}<\omega^\alpha\ln(\omega)$ $\ln(\ln(\omega))$

$<\exp(\omega^\alpha/\ln(\omega))$

$<\exp(\alpha\omega^\beta)/\omega^\gamma/\ln(\omega)^\delta<\exp(\omega)/\omega^\alpha/\ln(\ln(\omega))
<\exp(\alpha\omega^\beta)/\omega^\gamma<\exp(\omega)$ $\ln(\ln(\omega))/\omega^\alpha
<\exp(\alpha\omega^\beta)\ln(\omega)^\gamma/\omega^\delta<\exp(\omega)\ln(\ln(\omega))/\ln(\omega)
<\exp(\alpha\omega^\beta)/$ $\ln(\omega)^\gamma<\exp(\omega)\ln(\ln(\omega))/\ln(\omega)
<\exp(\alpha\omega^\beta)/\ln(\ln(\omega))<\exp(\alpha\omega^\beta)
<\exp(\alpha\omega^\beta)$ $\ln(\ln(\omega))<\exp(\omega)\ln(\omega)/\ln(\ln(\omega))
<\exp(\alpha\omega^\beta)\ln(\omega)^\gamma<\exp(\omega)\ln(\omega)\ln(\ln(\omega))$ $
<\exp(\alpha\omega^\beta)\omega^\gamma/\ln(\omega)^\delta<\exp(\omega)\omega^\alpha/\ln(\ln(\omega))
<\exp(\alpha\omega^\beta)\omega^\gamma<\exp(\omega)\omega^\alpha$ $\ln(\ln(\omega))
<\exp(\alpha\omega^\beta)\omega^\gamma\ln(\omega)^\delta$

$<\omega^{\omega^\alpha}<\exp(\exp(\omega)/\omega^\alpha)<\exp(\exp(\omega)/\ln(\omega))$

$<\exp(\exp(\omega)-\omega)/\omega^\alpha<\exp(\exp(\omega)-\omega)/\ln(\omega)
<\exp(\exp(\omega)-\alpha\omega^\beta)<\exp(\exp(\omega)-\omega)\ln(\omega)
<\exp(\exp(\omega)-\omega)\omega^\alpha$

$<\exp(\exp(\omega))/\omega^\alpha/\ln(\omega)<\exp(\exp(\omega))/\omega^\alpha
<\exp(\exp(\omega))\ln(\omega)/\omega^\alpha<\exp(\exp(\omega))/\ln(\omega)^\alpha
<\exp(\exp(\omega))/\ln(\ln(\omega))<\exp(\exp(\alpha\omega^\beta))
<\exp(\exp(\omega))$ $\ln(\ln(\omega))<\exp(\exp(\omega))\ln(\omega)^\alpha
<\exp(\exp(\omega))\omega^\alpha/\ln(\omega)<\exp(\exp(\omega))\omega^\alpha
<\exp(\exp(\omega))\omega^\alpha\ln(\omega)$

$<\exp(\exp(\omega)+\omega)/\omega^\alpha<\exp(\exp(\omega)+\omega)/\ln(\omega)
<\exp(\exp(\omega)+\alpha\omega^\beta)<\exp(\exp(\omega)+\omega)\ln(\omega)
<\exp(\exp(\omega)+\omega)\omega^\alpha$

$<\omega^{\exp(\omega)}<\exp(\exp(\omega)\omega^\alpha)<\exp(\exp(\exp(\omega)))$.

Also, not fitting neatly into the sequence above is: $\exp(\alpha\omega^\beta\pm\gamma\omega^\delta)/\omega^\epsilon
<\exp(\alpha\omega^\beta\pm\gamma\omega^\delta)/\ln(\omega)
<\exp(\alpha\omega^\beta\pm\gamma\omega^\delta)
<\exp(\alpha\omega^\beta\pm\gamma\omega^\delta)\ln(\omega)
<\exp(\alpha\omega^\beta$ $\pm\gamma\omega^\delta)\omega^\epsilon$ where $\beta>\delta$.\\

It seems likely that the entire set of Archimedean classes generated in this way can be ordered iteratively according to the following three rules.

\emph{Rule} 1. Can the class be sorted using the ordering of the real numbers $\alpha$, $\beta$, etc.? If not then go to Rule 2.

\emph{Rule} 2. Can the class be expressed in the form $fg$ or $f/g$? If so apply the $fg$ and $f/g$ rule. If not then go to Rule 3.

\emph{Rule} 3. Can the class be sorted using $\ln() < ()^\alpha$ ? If not then take the logarithm and return to Rule 1.\\

This iterative ordering can be turned into an equivalent deterministic ordering as follows:

\emph{Step} 1. Using the cardinal relation $f < \exp(f)$ construct the ordering of cardinals:
$\ldots< \ln\ln\ln(\omega) < \ln\ln(\omega) < \ln(\omega) < \omega < \exp(\omega) < \exp(\exp(\omega)) < \exp(\exp(\exp(\omega))) < \ldots$

\emph{Step} 2. Within each cardinal progressing from the innermost to the outermost exponential level (for $\omega$ and $\ln(\ldots)$ there is only one such level) sort by the $fg$ and $f/g$ rule:
 \begin{displaymath}
f^{\alpha_1}<f/g_1<f<fg_2<f^{\alpha_2}
\end{displaymath}
where $\alpha_1 < 1$, $\alpha_2 > 1$, sorting within categories by decreasing $g_1$ and by increasing $\alpha_1$, $g_2$ and $\alpha_2$.\\

Armed with these ordering relationships it is possible to show that what I've blithely called ``The infinite Archimedean classes" are in fact both infinite and Archimedean.

\emph{Lemma}: The ``infinite Archimedean classes" are infinite.

\emph{Proof}: The smallest prototype at generation $n$ is $\ln\ln…(\omega)$ with $n$ applications of ln. To show that this is smaller than the second smallest $\ln\ln\ldots(\omega)^\alpha$ with $n-1$ applications of ln, let $x = \ln\ln…(\omega)$ with $n-1$ applications of ln. From the previous lemma we already know that for sufficiently large $x$, $\ln(x) < x^\alpha$ for all $\alpha \in\mathbb{R}^+$. For rigorously proving that ``The smallest prototype at generation $n$ is $\ln\ln\ldots(\omega)$ with $n$ applications of ln", use induction from the previous generation. Now proof by contradiction. If $\ln\ln\ldots(\omega)$ with $n$ applications of ln is a real number $r$ then $\omega = \exp\exp\ldots(r)$ with $n$ applications of exp. But $\exp\exp\ldots(r)$ is a real number $< \omega$. Contradiction. So $\ln\ln\ldots(\omega)$ is infinite.

\emph{Lemma}: No two prototypes of the ``infinite Archimedean classes" belong to the same Archimedean class.

\emph{Proof}: The smallest prototype at generation $n$ is infinite. Since one way of creating a prototype at generation $n$ is by $f/g$ with $f > g$ from generation $n-1$, the ratio $f/g$ from generation $n-1$ must be larger than any real number. For any given natural number $m$ there is a larger real number, so at generation $n-1$ there does not exist an $m$ such that $f < mg$. So $f$ and $g$ cannot belong to the same Archimedean class. $n$ is arbitrary and finite, so no two prototypes can belong to the same Archimedean class.

\emph{Lemma}: If $f$ and $g$ are two prototypes of an infinite Archimedian class and  $\alpha\in\mathbb{R}^+$ then $f^{\alpha g}$ is also a prototype of an infinite Archimedian class.

\emph{Proof}: $f^{\alpha g}=\exp(g\ln(f))^{\alpha}$. Proof by induction. If $f$ and $g$ are of generation zero then $f^{\alpha g}=\exp(\omega\ln(\omega))^{\alpha}$ which is a prototype. If $f$ has the form $\omega$ or $\ln(\ldots)$ then it's immediately confirmed that $f^{\alpha g}$ is a prototype. If $f$ has the form $\exp(f_1)$ then $f^{\alpha g}=\exp(f_1g)$ is a prototype. If $f$ has the form $f_1^\beta$ for finite positive $\beta$ then $f^{\alpha g}=\exp(g\ln(f_1))^{\alpha\beta}$ so is is a prototype if all previous generations yield prototypes. If $f$ has the form $f_1f_2$ then $f^{\alpha g}=(\exp(g\ln(f_1))$ $\exp(g\ln(f_2)))^{\alpha}$ so is is a prototype if all previous generations yield prototypes. If $f$ has the form $f_1/f_2$ where $f_1 > f_2$ then $f^{\alpha g}=(\exp(g\ln(f_1))/\exp(g\ln(f_2)))^{\alpha}$ so it is a prototype if all previous generations yield prototypes. This exhausts all possibilities.\\

A third equivalent way to order the infinite Archimedean classes is to first order them as real numbers for sufficiently large $m$ and then take the limit of the ordering as $m$ tends to $\omega$.

\emph{Lemma}: Given generation $n$, maximum and minimum $\alpha_i$ and $|1-\alpha_i |$ for all real numbers $\alpha_i$ in the prototype of an infinite Archimedean class, there exists a finite $x$ such that for all $m > x$ the ordering of prototypes on $m$ matches that on $\omega$.

\emph{Proof}: Let $\Delta\alpha$ be the minimum $|1-\alpha_i |$. In order to preserve ordering, the inequalities that need to be satisfied are $x^{1-\Delta\alpha}<x/\ln(y)<x<x\ln(y)<x^{1+\Delta\alpha}$, $\ln(x)<x^{\alpha_{\min}}$, $x^{\alpha_{\max}}<exp(x)$ where $y$ depends on $x$ and generation $n$, and $\alpha$ represents any of the real numbers $\alpha_i$ in the prototype of an infinite Archimedean class. Each of these inequalities can be inverted to yield a value of $x$ exceeding the strict minimum necessary that depends on $n$, $\Delta\alpha$, $\alpha_{\min}$, $\alpha_{\max}$. Then for all $m > \max(x)$ the ordering of prototypes on $m$ matches that on $\omega$.\\

Returning now to the cardinal jump rule. This is important because it allows prototypes to reach beyond the end of the Aleph sequence of cardinals. The ordering is immediate. $\exp^{\omega}(f)$ is larger than any other manipulation of $f$ and $\ln^{\omega}(f)$ is between any other manipulation of $f$ and one.

\section{Limits of smooth sequences, and ordered fields}

What is a limit? What is the philosophy behind it? It seems to me that the limit of any sufficiently nice function $f(x)$ at the point $b$ must satisfy:
\begin{displaymath}
\lim_{x\to b}f(x)=f(b)
\end{displaymath} 
and in particular
\begin{displaymath}
\lim_{x\to b}x=b
\end{displaymath}
no matter what number system the function and values reside on. The difficulty lies in defining exactly what a `sufficiently nice function' is.

The normal (Cauchy) limit can be defined as~\cite{20}:
\begin{itemize}
\item	An infinite sequence of real numbers is a map $s:\mathbb{N}\to\mathbb{R}$ and denote $s(n)$ by $s_n$. This sequence converges to the real number $c$ if given any real number $\epsilon > 0$, there exists an integer $N_\epsilon$ such that $|s_n-c|<\epsilon$ for all $n\ge N_\epsilon$.
\end{itemize}
The new type of limit introduced here is a form of “asymptotic limit” where the asymptotes are prototypes of Archimedean classes. Let $p$ denote the prototype of Archimedean class $O(p_n)$. Define the \textbf{leading term limit} by:
\begin{itemize}
\item	If for every positive real number $\epsilon$, there is an integer $N_\epsilon$, a non-zero real number $c$ and an Archimedean class $p$, such that $|(s_n/p_n)-c|<\epsilon$ for all integer $n\ge N_\epsilon$, then it is said that the leading term limit $^1\lim s_n=cp$. If, instead, $s_n=0$ for all $n \ge N_\epsilon$ then $^1\lim s_n=0$.
\end{itemize}
That gives the leading term of the limit, and that will suffice for many purposes. To get the second term subtract off the leading term and continue the process.
\begin{itemize}
\item	If  $^1\lim s_n=cp$ and for every positive real number $\epsilon$, there is an integer $N_\epsilon$, a non-zero real number $d$ and Archimedean class $q$, such that $|((s_n-cp_n)/q_n)-d|<\epsilon$ for all integer $n \ge N_\epsilon$, then it is said that the \textbf{second term limit}  $^2\lim s_n=dq$.
\end{itemize}
For the third term of the limit subtract off the second term as well and repeat.
\begin{itemize}
\item $\lim s_n={}^1\lim s_n+{}^2\lim s_n +\ldots$ as far as you want to take it.
\end{itemize}
Before going further, I want to explain what ``as far as you want to take it" means. Recall that L'H\^{o}pital's rule~\cite{21} uses derivatives to evaluate expressions that otherwise reduce to 0/0. If the first derivative fails then you want to take it to the second derivative, etc. Similarly, one of the applications of infinite numbers is to evaluate expressions that otherwise evaluate to $\infty-\infty$. If the leading term fails then you want to take it to the second term, etc. That's what ``as far as you want to take it" means.

\emph{Lemma}: The leading term limit agrees with the Cauchy limit when the Cauchy limit is nonzero.

\emph{Proof}: $p = p_n = 1$ and the leading term limit reduces exactly to the Cauchy limit.\\

Let $^1c,{}^2c,{}^1p,{}^2p$ denote $c, d, p, q$ and let $I$ be a finite natural number. Then the \textbf{limit} is:
\begin{displaymath}
\lim_{n\to\omega}s_n=\sum_{i=0}^I{}^ic{}^ip\quad \mbox{ or }0
\end{displaymath}

There have been many approaches to infinite and infinitesimal numbers in the past that  parallel the present approach, using a series of products of a finite number and a prototype. Two of the best known are the Hahn series~\cite{9,22,Hahn} and the Levi-Civita series~\cite{18,23}, both are subsets of the hyperreals.  Both use power series and allow the summation to continue to infinity. The reason for insisting here that the sum be finite is that a finite sequence of $N_\epsilon$ must have a finite upper bound. This upper bound forces the numbers generated here to be either infinite, infinitesimal or finite real. Because both Hahn series and Levi-Civita series are infinite, they allow the existence of numbers that do not fall into these three categories, numbers such as $\sin(\omega)$ that are finite but not real. The new definition of limit presented in this paper deliberately excludes all such numbers, although a different limit could be defined to include them (think ``limit cycle").

Using the definition of limit above not every sequence will have even a leading term limit (I'll show a way around this problem in Part 2 of this paper). A name is needed for those sequences that do, perhaps \textbf{leading limitable sequences}. Leading limitable sequences do not form a field because they are not closed under term by term addition. Both $\exp(n) + \sin(n)$ and $-\exp(n)$ are leading limitable sequences but their sum $\sin(n)$ is not. Use the word \textbf{smooth} for those sequences that have limits using the definition above out as far as you want to take it. These also do not form a field. Use the word \textbf{oscillatory} for all sequences that are not smooth.\\

In what follows, unless otherwise stated, $n\to\omega$ is implied. And denote the limit of a sequence $\lim_{n\to\omega}s_n$ by $\lim s$.

Denote the \textbf{infinite numbers} by $\mathbb{I}_n$ and let all the numbers in it be constructed using
\begin{displaymath}
\lim s/\lim t
\end{displaymath}
where $s$ and $t$ are smooth sequences and $t$ has a nonzero limit. The symbol $\mathbb{I}_n$ is chosen becuse it stands for all three of Infinite, Infinitesimal, and sequence. It is a temporary notation until this set of numbers can be absorbed into pre-existing descriptions of the infinite.

\emph{Lemma} Every element of $\mathbb{I}_n$ is a limit of a sequence.

\emph{Proof} Every element of $\mathbb{I}_n$ has form $\lim s / \lim t$ which is the limit of the sequence $\lim_{n\to\omega}\lim s_n / \lim t_n$\\

Is $\mathbb{I}_n$ complete? Not without the feedback rule defined below. Without this rule it is found that $\exp(\exp(\omega)/(\omega+1))$ ought to be included among the prototypes but is not yet there.

Add this rule to the definition of infinite prototypes:
\begin{itemize}
\item If $\lim t = \sum {}^id\:{}^iq$ is any positive limit with no infinitesimal terms and $f$ is any infinite prototype that satisfies $\ln(f) > (\lim t)^n$ for every integer $n$ then $f^{1/(\lim t)}$ is an infinite prototype. Call this the \textbf{feedback rule}.
\end{itemize}

\emph{Lemma} When $g\in\mathbb{I}_n$ and the conditions for the feedback rule are satisfied, $f^g$ is a prototype (infinite, 1 or infinitesimal).

\emph{Proof} $g=0$ in which case $f^g=1$, or $g=({}^1c{}^1p+{}^2c{}^2p+\ldots{}^mc{}^mp)/({}^1d{}^1q+{}^2d\:{}^2q+\ldots{}^nd{}^nq)$. Rewrite this as $g=\pm({}^mp/{}^nq)({}^1c'{}^1p'+{}^2c'{}^2p'+\ldots{}^mc')/({}^1d'{}^1q'+{}^2d'\:{}^2q'+\ldots{}^nd')$ where ${}^ic'={}^ic\left|{}^1c\right|/{}^1c$ and ${}^jd'={}^jd\left|{}^1d\right|/{}^1d$. The top and bottom limits are now each positive and have no infinitesimal terms.
From a previous lemma, $f^{cp}$ is an infinite prototype for all infinite prototypes $f$ and $p$ when $c\in\mathbb{R}^+$. If $p=1$ then this is also true by the definition of infinite prototype. If $c<0$ then $f^{cp}=1/f^{-cp}$ is an infinitesimal prototype. Because the product of two prototypes is a prototype, the finite product $\prod_i f^{\:{}^ic'\:{}^ip'}$ is a prototype. Let $f'=\prod_i f^{\:{}^ic'\:{}^ip'}$ and $\lim t={}^1d'{}^1q'+{}^2d'\:{}^2q'+\ldots{}^nd'$ which is positive with no infinitesimal terms. So when $\ln(f') > (\lim t)^n$ for every integer $n$ the feedback rule says that ${f'}^{1/(\lim t)}$ is a prototype. Nearly there.
$f^g=({f'}^{1/(\lim t)})^{(\pm{}^mp/{}^nq)}$. Prototypes are closed under division and under the action $p^{-1}$. If ${}^mp/{}^nq$ is infinite then the definition of infinite and infinitesimal prototypes says that the result is a prototype of the same type. If ${}^mp/{}^nq$ is infinitesimal then the result is either 1 (eg. $\omega^{1/\omega}$ yields the prototype 1) or a prototype of the same type.\\

With the feedback rule, closure under multiplication and division still follow automatically because the rules allowing $fg$ and $f/g$ to be generated are unchanged.

Ordering of prototypes with the feedback rule follows from the expansion of $1/(\lim t)$ as a power series. $\lim t={}^1d{}^1q+{}^2d\:{}^2q+\ldots{}^nd{}^nq$. Write this as $\lim t={}^1d{}^1q(1+{}^2d'\:{}^2q'+\ldots{}^nd'{}^nq)={}^1d{}^1q(1-x)$ where $|x|<1$. Then $1/(\lim t)=(1/{}^1d{}^1q)(1+x+x^2+x^3+\ldots)$. Any finite part of the infinite series equalling $f^{1/(\lim t)}$ can be reproduced by successive uses of the operations $f^g$, $f^\alpha$, $fg$ and $f/g$ to build it up term by term. Let's say that finite operations build up to $(1+x+x^2+\ldots+x^n)$ where  $n$ is odd. Then the next term $x^{(n+1)}>0$ and so $f^{(1/{}^1d{}^1q)(1+x+x^2+\ldots+x^n)}<f^{1/(\lim t)}<f^{(1/{}^1d{}^1q)(1+x+x^2+\ldots+(1\pm\epsilon)x^n)}$ for arbitrarily small $\epsilon$ with the sign depending on the sign of $x^n$. That's the ordering relationship. Because $f^{1/(\lim t)}$ always fits between two existing infinite classes it immediately follows that it is both infinite and positive.\\

\emph{Lemma} $f^{1/(\lim t)}$ builds a new Archimedean class.

\emph{Proof} This is true if $r=f^{1/(\lim t)}/f^{(1/{}^1d{}^1q)(1+x+x^2+\ldots+x^n)}>m$ for all $m\in\mathbb{N}$ because the ratio on the other side is larger because it relates to an earlier term of the power series. $\ln(r)/\ln(f)=1/(\lim t)-(1/{}^1d{}^1q)(1+x+x^2+\ldots+x^n)=(1/{}^1d{}^1q)(x^{(n+1)}+x^{(n+2)}+\ldots)=x^{n+1}/\lim t$. Both $x^{n+1}=|x|^{n+1}$ and $\lim t$ are positive by definition. So this construction builds a new Archimedean class if $x^{n+1}\ln(f)/\lim t>m$ for all  $m\in\mathbb{N}$. From the definition of the feedback rule, $\lim t$ has no infinitesimal terms, so $|x|>1/\lim t$. So $x^{n+1}\ln(f)/\lim t>\ln(f)/(\lim t)^{n+2}$. From the definition of the feedback rule, $\ln(f)>(\lim t)^k$ for all $k$. Set $k=n+3$. Then $\ln(f)/(\lim t)^{n+2}>\lim t$ is greater than any natural number. So  $f^{1/(\lim t)}$ builds a new Archimedean class.\\

The following lemmas clarify the status of limits and lead to the conclusion that $\mathbb{I}_n$ is an ordered field.

\emph{Lemma}: If elements of a nonzero Archimedean class are defined using the above definition of limit then that Archimedean class is isomorphic to the nonzero reals.

\emph{Proof}: For every nonzero real number $c$ a sequence $s_n=cp_n$ has leading term limit $cp$. So the Archimedean class contains a subset isomorphic to the nonzero reals. Hilbert has proved that the reals are isomorphic to the largest Archimedean class. So the Archimedean class is isomorphic to the nonzero reals.

\emph{Corollary}: Elements of nonzero Archimedean classes defined using the above definition of limit are identical to those defined using leading term limits.

\emph{Lemma}: If elements of a nonzero Archimedean class are defined by the limit of a sequence of rationals then that Archimedean class is isomorphic to the nonzero reals.

\emph{Proof}: Let $c_n$ be a Cauchy sequence of rationals converging to nonzero real number $c$ and let $s_n=c_n p_n$. Then the leading term limit $|(s_n/p_n)-c|=|c_n-c|<\epsilon$ for all integer $n\ge N_\epsilon$ and so the leading term limit is $cp$. The limit of Cauchy sequences of rationals defines the real numbers, so the Archimedean class is isomorphic to the nonzero reals.

\emph{Lemma}: The set of all numbers defined by leading term limits is ordered.

\emph{Proof}: The Archimedean classes defined using commutative Big O notation are strictly ordered and the real numbers are ordered within each Archimedean class. If two numbers defined by leading term limits are in different Archimedean classes then the ordering of Archimedean classes is used. If they fall in the same Archimedean class then the ordering of the nonzero reals is used.

\emph{Lemma}: The set of all limits is ordered.

\emph{Proof}: Order first by leading term limits. If they are equal, order by second term limit. If they are equal, order by third term limit. As far as you want to take it.

\emph{Lemma}: Nonzero leading term limits are closed under multiplication and division.

\emph{Proof}: The set of Archimedean classes excluding the singleton class of zero is closed under multiplication and division. The nonzero real numbers are closed under multiplication and division. The product and quotient of two leading term limits is $(cp)(dq)=(cd)(pq)$ and $(cp)/(dq)=(c/d)(p/q)$ which are also leading term limits.

\emph{Lemma}: Limits of smooth sequences possess all the properties of a field other than a multiplicative inverse.

\emph{Proof}: The proof is elementary, let's take the distributive law as an example. Denote prototypes by $p$, $q$ and $r$. Denote real numbers by $c$, $d$ and $e$. Denote sequences by $s$, $t$ and $u$. Let $\lim s = c_1p_1 + \ldots + c_i p_i$ and $\lim t = d_1 q_1 + \ldots + d_j q_j$ and $\lim u = e_1 r_1 + \ldots + e_k r_k$.

$\lim s \times (\lim t + \lim u)
= c_1p_1(d_1 q_1 + \ldots + d_j q_j + e_1 r_1 + \ldots + e_k r_k) + \ldots + c_i p_i (d_1 q_1 + \ldots + d_j q_j + e_1 r_1 + \ldots + e_k r_k)
= c_1p_1(d_1 q_1 + \ldots + d_j q_j) + c_1p_1(e_1 r_1 + \ldots + e_k r_k) + \ldots + c_i p_i (d_1 q_1 + \ldots + d_j q_j) + c_i p_i (e_1 r_1 + \ldots + e_k r_k)
= c_1p_1(d_1 q_1 + \ldots + d_j q_j) + \ldots + c_i p_i (d_1 q_1 + \ldots + d_j q_j) + c_1p_1(e_1 r_1 + \ldots + e_k r_k) + \ldots + c_i p_i (e_1 r_1 + \ldots + e_k r_k)
= \lim s \times \lim t + \lim s \times \lim u$.

The multiplicative identity has $c_1 = p_1 = 1$ and $i = 1$, and this is unique because $p_1$ only has the options of being 0, 1, infinite or infinitesimal, and $c_1$ must be finite. The multiplicative inverse is not present when $i > 1$.

\emph{Lemma}: $\mathbb{I}_n$ is ordered.

\emph{Proof}: $\lim s$ is closed under multiplication and ordered. So $\lim s \lim t$ is ordered. $\lim t$ is nonzero. Without loss of generality take $\lim t > 0$. Then $(\lim s/\lim t)_1> (\lim s/\lim t)_2$ if and only if $(\lim s)_1(\lim t)_2>(\lim s)_2(\lim t)_1$.

\emph{Lemma}: $\mathbb{I}_n$ is a field.

\emph{Proof}: The proof is elementary, again let's take the distributive law as an example.
\begin{eqnarray*}
&(\lim s/\lim t)_1\times((\lim s/\lim t)_2+(\lim s/\lim t)_3)\\
&=(\lim s/\lim t)_1\times((\lim s)_2(\lim t)_3+(\lim s)_3(\lim t)_2)/((\lim t)_2(\lim t)_3)\\
&=((\lim s)_1(\lim s)_2(\lim t)_3+(\lim s)_1(\lim s)_3(\lim t)_2)/((\lim t)_1(\lim t)_2(\lim t)_3)\\
&=((\lim s)_1(\lim s)_2)/((\lim t)_1(\lim t)_2)+(\lim s)_1(\lim s)_3)/((\lim t)_1(\lim t)_3)\\
&=(\lim s/\lim t)_1\times(\lim s/\lim t)_2+(\lim s/\lim t)_1\times(\lim s/\lim t)_3
\end{eqnarray*}

\setlength{\tabcolsep}{3 pt}
\begin{table}
\caption{Summary of Field Properties.}
\vspace{10pt}
\begin{tabular}{|l|cccc|}
\hline
& Prototype & Leading lim & Limit & $\mathbb{I}_n$ \\
& $p$ & $^1\lim = cp$ & $\lim\! =\! \Sigma cp$	 & $\lim s\! /\! \lim t$ \\\hline
Closure under addition	&$\times$&$\times$&\checkmark&\checkmark\\
Closure under multiplication&\checkmark&\checkmark&\checkmark&\checkmark\\
Associativity of addition&&&\checkmark&\checkmark\\
Associativity of multiplication	&\checkmark&\checkmark&\checkmark&\checkmark\\
Commutativity of addition	&&&\checkmark&\checkmark\\
\begin{small}Commutativity of multiplication\end{small}&\checkmark&\checkmark&\checkmark&\checkmark\\
Additive identity&&&\checkmark&\checkmark\\
Multiplicative identity&\checkmark&\checkmark&\checkmark&\checkmark\\
Additive inverse&&&\checkmark&\checkmark\\
Multiplicative inverse&\checkmark&\checkmark&$\times$&\checkmark\\
Distributivity&&&\checkmark&\checkmark\\
Ordered&\checkmark&\checkmark&\checkmark&\checkmark\\\hline
\end{tabular}
``prototype" denotes the set of prototypes of all Archimedean classes.\\
\end{table}

\section{Analogies between the infinites and the surreals}

A very brief summary of the relevant properties of the surreals is as follows. A surreal number $x\in \bf{N_O}$ can be written in the form $\{X_L|X_R\}$ where $X_L$ and $X_R$ are sets of surreal numbers. Let $x_L\in X_L$ and $x_R\in X_R$. All $x_L\le x_R$. The reals are a subset of the surreals, the product of two surreal numbers is a surreal number. The division of two surreal numbers, say $z=x/y$ can be constructed using $z=\{Z_L|Z_R\}$ where $\forall z_L\in Z_L:yz_L<x$ and $\forall z_R\in Z_R:yz_R>x$. The successor of $x$ is $x+1=\{x|\}$ where the null set is written as a blank. The surreals are ordered in generations, and if the expression $\{X_L|X_R\}$ allows one or more possibilities from an earlier generation then it is set to the value from the earliest generation. The earliest generation is $0=\{|\}$. The earliest generation infinite numbers are $\omega=\{\mathbb{Z}|\}$ and $-\omega=\{|\mathbb{Z}\}$.\\

It's instructive to compare the present momenclature with that of Alling (1987) ``Foundations of analysis over surreal number fields".\cite{7}

\setlength{\tabcolsep}{3 pt}
\begin{table}[ht]
\caption{Rough matching of present nomenclature and that of Alling(1987)}
\vspace{10pt}
\centering
\begin{tabular}{|l|l|}
\hline
Alling(1987)&This paper\\\hline
Pseudo-limit&Limit\\
Pseudo-convergent&Smooth\\
Limit&Leading term limit\\
U&$\mathbb{R}\backslash 0$\\
Value group G&$\mathbb{R}$\\
K $=\mathbf{N_O}=$ Surreals&$\mathbb{I}_n$\\
Valuation V&Mapping from $\mathbb{I}_n$ to $\mathbb{R}\backslash 0$\\
Ordinal number $\mathbf{O_n}$&$\mathbb{N}$\\
Valuation ring \bf{O}&Finites + Infinitesimals\\
A $=\sum_{\alpha<\lambda}a_{\alpha}\omega^{g(\alpha)}$&$\lim s=\sum_{i=0}^I{}^ic{}^ip$\\\hline
\end{tabular}
\end{table}

Alling ascribes the invention of the pseudo-limit on psuedo-convergent sequences to Ostrowski(1935).\cite{Ostrowski} His definition of ``limit" as ``the unique simplest pseudo-limit" may not correspond all that closely to my ``leading term limit". The surreals $\mathbf{N_O}$ may not match the infinites $\mathbb{I}_n$, that has yet to be determined and is discussed below. Alling's use of valuation V serves the same purpose as the present prototypes, but in reverse. A valuation is a mapping from the non-Archimedean field K to the nonzero reals; a prototype is a mapping from the nonzero reals to a set of values on the non-Archimedean field K. The ordinal numbers $\mathbf{O_n}$ are vaster than the natural numbers $\mathbb{N}$, but I've matched them in the table because Alling uses $\mathbf{O_n}$ the way this paper uses $\mathbb{N}$. After this, Alling goes on to define analytic functions of a surreal variable, such as trigonometric fuctions.\\

Direct analogies that exist between $\mathbb{I}_n$ and $\mathbf{N_O}$ include the following.

\emph{Analogy} The infinite Archimedean classes contain $\omega$ and are closed under a finite number of operations of  $\ln( )$, $\exp( )$, finite positive power, a multiplication of two unequal classes, and a division of a larger class by a smaller.

\emph{Lemma} The infinite surreals contain $\omega$ and are closed under a finite number of operations of  $\ln( )$, $\exp( )$, finite positive power, a multiplication of two unequal surreals, and a division of a larger surreal by a smaller.

\emph{Proof} The infinite surreals are closed under pairwise multiplication and division and contain $\omega$. Let $x$ be an infinite surreal number and define $\exp(x)=\{1,1+x,1+x+x^2/2,1+x+x^2/2+x^3/6,\ldots|\}-1$ The ``$-1$" is to remove the ``$+1$" added by the successor function. Because surreal numbers are closed under pairwise addition and pairwise multiplication, any finite power series of a surreal number is a surreal number. Infinite sets are allowed so $\exp(x)$ is a valid surreal number. Define $y=\ln(x)=\{Y_L|Y_R\}$ where $\forall y_L\in Y_L:\exp(y_L)<x$ and $\forall y_R\in Y_R:\exp(y_R)>x$. Then $\ln(x)$ is a valid surreal number. Define $x^\alpha=\exp(\alpha\ln(x))$. $x^\alpha$ is a valid surreal number because the surreal numbers are closed under $\exp$, $\ln$ and pairwise multiplication. Both $\exp(x)$ and $\ln(x)$ are infinite.

\emph{Analogy} The infinite Archimedean classes are closed under operations $\exp\exp\exp\ldots() =$ $ \exp^{\omega}()$ and $\ln\ln\ln\ldots() =$ $ \ln^{\omega}()$.

\emph{Lemma} The infinite surreals are closed under operations $\exp\exp\exp\ldots()$ and $\ln\ln\ln\ldots()$.

\emph{Proof} Let $x$ be an infinite surreal and define $\exp\exp\exp\ldots(x)=\{x,\exp(x),$ $\exp\exp(x),\exp\exp\exp(x),\ldots|\}-1$ and $\ln\ln\ln\ldots(x)=\{x,\ln(x),\ln\ln(x),$ $\ln\ln\ln(x),\ldots|\}-1$. Each term of the sequence exists by closure under exp and ln, and so the limit exists as an infinite surreal number.

\emph{Analogy} Prototypes include 0, 1, the reals and 1/the-infinite-Archimedean-classes.

\emph{Lemma} Surreals include 0, 1, the reals and $1/x$ whenever $x$ is a surreal number.

\emph{Proof} The surreals include 0, 1 and the reals. Define $y=1/x=\{Y_L|Y_R\}$ where $\forall y_L\in Y_L:xy_L<1$ and $\forall y_R\in Y_R:xy_R>1$. Then $1/x$ is a valid surreal number.

\emph{Analogy} Let ${}^iq$ be a prototype and ${}^id$ be real. $\lim t = \sum {}^id\:{}^iq$ is any positive limit with no infinitesimal terms and $f$ is any infinite prototype that satisfies $\ln(f) > (\lim t)^n$ for every integer $n$ then $f^{1/(\lim t)}$ is an infinite prototype.

\emph{Lemma} For surreals ${}^iq$ and ${}^id$ and $\lim t = \sum {}^id\:{}^iq$, the value $f^{1/(\lim t)}$ is a surreal number.

\emph{Proof} The surreals are closed under pairwise addition and multiplication so $\lim t$ is a surreal. $f^{1/(\lim t)}=\exp(\ln(f)/\lim t)$ and because of closure under exp, ln and pairwise division this is a surreal number.

\emph{Analogy} denote the infinite numbers by $\mathbb{I}_n$ and let all the numbers in it be constructed using $\lim s/\lim t$

\emph{Lemma} $\lim s/\lim t$ is a valid surreal number.

\emph{Proof} This follows immediately from closure under pairwise division.\\

By this stage the reader may be thinking "So what, surreal numbers can do everything". If so, take note of the following.

\emph{Analogy} The infinite numbers exclude numbers like $\sin(\omega)$ that are finite but not real.

\emph{Lemma} $\sin(\omega)$ is not a separate surreal number because $\sin(\omega)=0$. Numbers derived from $\sin(\omega)$ by many mathematical operations are also excluded from the surreals.

\emph{Proof} $\sin(x)$ can be expressed as a power series $\sin(x)=x-x^3/3!+x^5/5!-x^7/7!+\ldots$. It is well known that every partial sum changes from being greater to $\sin(x)$ to less than $\sin(x)$ or vice versa so we can define $\sin(x)$ on the surreals using $y=\sin(x)=\{x-x^3/3!,x-x^3/3!+x^5/5!-x^7/7!,\ldots|x,x-x^3/3!+x^5/5!,\ldots\}$ for $x>0$ and $y=\sin(x)=\{x,x-x^3/3!+x^5/5!,\ldots|x-x^3/3!,x-x^3/3!+x^5/5!-x^7/7!,\ldots\}$ for $x<0$. When $x$ is infinite, every term $y_L\in Y_L<0$ and every term $y_R\in Y_R>0$ so $\sin(x)$ can equal zero. The value is allocated to the earliest generation consistent with the definition, so $\sin(x)=0$ when $x$ is infinite. Values derived from the definition of $\sin(x)$ by finite addition, multiplication and division preserve the split between $Y_L$ and $Y_R$ and so are the same as the equivalent operations on 0.\\

All the above suggests that $\mathbb{I}_n\subset\bf{N_O}$. Ideally the two are equivalent, and if $\mathbb{I}_n$ falls short then extra rules, possibly of the feedback type, should be constructed to fill the gap. Two properties of the surreals that might be portable to the infinites are that on the surreals every infinite sequence has well-defined sup and sub values, and many functions on the surreals have well-defined inverse functions.

It will be beiefly proved in Part 2 of this paper that $\mathbb{I}_n$ is a proper subset of the hyperreals $^*\mathbb{R}$. The two are not equivalent. $\mathbb{I}_n$ and $\bf{N_O}$ are both constructed without the need for the axiom of choice. $^*\mathbb{R}$ cannot be defined without the axiom of choice.

\section{Conclusions}

Part 1 of this paper introduces infinite (and infinitesimal) numbers as a collection of Archimedean classes, and defines each class as the asymptotic limit of a specific prototype sequence. This approach seems to be new, and highly practical because infinite sequences appear in many places in pure and applied mathematics. In Part 1 only a relatively few sequences, those that are called ``smooth" or are the ratio of two smooth sequences, are allocated limits. The resulting infinite numbers, called here $\mathbb{I}_n$, are shown to be an ordered field under the normal operations of addition and multipication.

This construction can be read in concert with earlier constructions of infinite and infinitesimal numbers: the hyperreals, the surreals, by consideration of the geometric continuum, and from power series. The one that most closely resembles the present construction is the surreals.

The paper leaves some questions unanswered. Proper proofs/disproofs are needed for each of the following conjectures:
\begin{itemize}
\item	The infinite numbers presented here are a subset of the surreal numbers.
\item	The ordering relationship for prototypes presented above is correct.
\item	$\mathbb{I}_n$ can be completely defined by smooth sequences of rationals.
\end{itemize}
Other pure mathematicians may wish to consider extensions involving topology, field theory, sequences on hypernaturals, differential vector calculus, etc.\\

The use of smooth sequences in the construction of $\mathbb{I}_n$ provides the springboard for Part 2 of this paper, on the limits of oscillatory sequences and on some elementary applied mathematics applications of limits on sequences.

\end{document}